\documentclass[letterpaper, 10 pt, conference]{ieeeconf}

\IEEEoverridecommandlockouts
\usepackage{subfig}
\usepackage{float}
\usepackage{transparent}
\usepackage{amsmath}
\usepackage{amssymb}
\usepackage{graphicx}
\usepackage{color}
\usepackage{balance}

\usepackage{mathtools}

\mathtoolsset{showonlyrefs}

\def\rT{{\rm T}}        

\def\<{\leqslant}           
\def\>{\geqslant}           

\def\[{[\![}
\def\]{]\!]}

\usepackage{transparent}
\usepackage{amsmath}
\usepackage{amssymb}
\usepackage{graphicx}
\usepackage{color}
\usepackage{fullpage}

\newcommand{\cI}{{\cal I}}
\newcommand{\cJ}{{\cal J}}

\newcommand{\cL}{{\cal L}}
\newcommand{\cM}{{\cal M}}

\newcommand{\gamovertwo} {{\frac{\gamma^2}{2}}}
\newcommand{\gammovertwo} \gamovertwo
\newcommand{\gammuovertwo} \gamovertwo
\newcommand{\gammuepstovertwo} \gamovertwo
\newcommand{\gammubarepstovertwo} \gamovertwo

\newcommand{\noncr} {\nonumber\\}

\newcommand{\beasnum}{\begin{eqnarray}}
\newcommand{\eeasnum}{\end{eqnarray}}
\newcommand{\beas}{\begin{eqnarray*}}
\newcommand{\eeas}{\end{eqnarray*}}

\newcommand{\be}{\begin{equation}}
\newcommand{\ee}{\end{equation}}
\newcommand{\ba}{\begin{array}}
\newcommand{\ea}{\end{array}}

\newtheorem{theorem}            {Theorem}[section]

\newtheorem{sideremark}         [theorem]{Remark}
\newtheorem{sideeg}           [theorem]{Example}
\newtheorem{sideconj}           [theorem]{Conjecture}

\def\argmin                     {\mathop{\rm argmin}}

\newcommand{\qed} {\hskip 0.2em\lower 0.7ex\hbox{\vbox{\hrule
\hbox{\vrule height 1.2ex\hskip 0.4em\vrule height 1.2ex}
\hrule}}}

\newcommand {\qform}[2]{{#1}^T {#2} {#1}}
\newcommand {\qformvec}[2]{\left(\begin{array}{cc}{#1}^T & 1\end{array}\right)  {#2} \left(\begin{array}{c}{#1} \\1\end{array}\right)}

\newcommand {\braket}[2]{\langle{#1}, {#2}\rangle}

\def \qweight {Q_\eta}

\def \minplus {\bigoplus}

\usepackage{color}

\def \figname {fig:minplusrobustest:}
\def \eqnname {eq:minplusrobustest:}

\begin{document}

\title{
Min-Plus Techniques for Set-Valued State Estimation
}
%

\author{Abhijit G. Kallapur \thanks{S.~Sridharan and W.M.~McEneaney are with the Department of Mechanical and Aerospace Engineering,
University of California, San Diego, CA 92093-0411, USA. Their research is supported under AFOSR grant FA9550-10-1-0233.
        { { srsridharan@eng.ucsd.edu, wmceneaney@eng.ucsd.edu}}}%
\and
 Srinivas Sridharan
\thanks{A. Kallapur and I. Petersen are with the School of Engineering and Information Technology,
University of New South Wales at the Australian Defence Force Academy, Canberra, ACT 2600, Australia         { a.kallapur@adfa.edu.au, i.petersen@adfa.edu.au}.  Their  work was supported by the Australian Research Council} 
          \and
 William M. McEneaney
\and 
 Ian R. Petersen
}

\maketitle

\begin{abstract}
This article approaches deterministic filtering via an application of the  min-plus linearity of the corresponding dynamic programming operator.  This filter design method  yields  a set-valued state estimator for discrete-time nonlinear systems (nonlinear dynamics and output functions). The energy bounds in the process and the measurement disturbances are modeled using a sum quadratic constraint. The filtering problem is recast into an optimal control problem in the form of a Hamilton-Jacobi-Bellman (HJB) equation, the solution to which is obtained by employing the min-plus linearity  property of the dynamic programming operator. This approach enables the solution to the HJB equation and the design of the filter without recourse to linearization of the system dynamics/ output equation.
\end{abstract}

\section{Introduction}
Deterministic filtering methods have been developed in the literature for linear and nonlinear systems as an alternative to stochastic techniques. They are especially applicable to situations where the noise characteristics are not stochastic and/or whose statistics are not known apriori. In such cases, the noise is typically modeled as an unknown process satisfying some bound in an $H_2$ or $H^{\infty}$ sense,
where this may be interpreted as a generalization of an energy bound; e.g., see \cite{mortensen1968,bertsekas1971,baras1995,fleming1997,mceneaney1998,james1998}. In particular, the set membership state estimation approach in \cite{bertsekas1971} provides a deterministic interpretation of the Kalman filter in terms of a set-valued state estimate, where the solution to the estimation problem is obtained by constructing the set of all possible states consistent with the given measurements. This set membership approach has been extended to estimation for nonlinear systems as presented in \cite{james1998,scholte2003,zhou2008,yang2009,ieeetac2009,acc2009}.

Most estimation schemes proposed for nonlinear systems both in the stochastic and the deterministic settings, use some kind of approximation schemes for the state dynamics which often consist of linearizing the state and measurement equations about a suitable operating point. Indeed, this is not an issue for systems with small nonlinearities but the effect of nonlinearity induced errors needs to be considered for systems with large nonlinearities as presented in \cite{scholte2003}. Another approach to nonlinear filtering that does not consider linearization of the underlying nonlinear dynamics is presented in \cite{fleming2000} using the max-plus machinery. There, the nonlinear filtering problem is recast into an optimal control problem  leading to  a Hamilton-Jacobi-Bellman (HJB) equation. 
We note that  the approach taken in \cite{fleming2000} bears a resemblance to that here in that max-plus machinery was employed
for value function propagation.
In both cases,  the value, as a function of the state variable is semiconvex,
where the space of semiconvex functions is a max-plus linear space (or moduloid),  for which  a countable set of quadratic functions
forms a basis (or spanning set).
However, in the case of \cite{fleming2000}, the basis elements used in the max-plus value function expansion were fixed.
Here, we adopt a modification of the more recently developed curse-of-dimensionality-free approach,
where that approach to infinite time-horizon control is adapted to the time-dependent filter value function propagation.
That approach has been demonstrated to be highly effective from a computational standpoint in several applications \cite{mceneaney2008curse}.
In particular, the quadratic functions used in the truncated max-plus expansion are boot-strapped by the algorithm.
We note that, in order to maintain computational tractability with this approach,
one employs a max-plus optimal projection at each time-step.
This is optimally performed by pruning the set of quadratics in the representation
\cite{mceneaney2009complexity}.

In this paper, we present a set-valued state estimation approach to nonlinear filtering for systems  with nonlinear dynamics and observations using min-plus methods to obtain the corresponding deterministic filter. The constraint on the system noise is described by a sum quadratic constraint (SQC) \cite{ieeetac2009,acc2009}. A set-valued state estimation scheme is utilized to reduce the filtering problem to a corresponding optimal control problem in terms of an HJB equation. The optimization problem consists of computing the minimum quadratic supply needed to drive the system to a given terminal state subject to the SQC. The computations are achieved by applying a min-plus scheme to the optimization process where the solution operator is linear in the min-plus algebra. Indeed, this scheme does not employ the linearization of the system dynamics and provides a less conservative solution in terms of the filter recursion equations.

The rest of the paper is organized as follows: Section \ref{sec:problem-formulation} describes the formulation of a nonlinear system with the noise bounded by an SQC. Section \ref{sec:svse} describes the set-valued state estimation scheme for nonlinear filtering, and recasts the nonlinear filtering problem into a corresponding optimal control problem.  The solution to the optimal control problem using the min-plus linearity property and the corresponding filter recursion equations which arise therefrom are discussed in Section \ref{sec:filter-equations}. An illustrative example is presented in Section \ref{sec:ex}, and the paper is concluded with remarks on future research in Section \ref{sec:conclusions}.

\section{Problem Formulation}
\label{sec:problem-formulation}
Consider a continuous-time system described by 
\begin{align}
\label{eq:ct_x}
	\dot{x}(t) &= a_c(t,x(t),u(t)) + D_c(t) \: w(t); \\
\label{eq:ct_y}
	y(t) &= c_c(t,x(t)) + v(t);
\end{align}
where $x(\cdot)\in\mathbb{R}^n$ is the state, $u(\cdot) \in \mathbb{R}^m$ is the known control input, $w(\cdot)\in\mathbb{R}^p$ and $v(\cdot)\in\mathbb{R}^l$ are the process and measurement disturbance inputs respectively, and $y(\cdot)\in\mathbb{R}^l$ is the measured output. $a_c(\cdot)$ and $c_c(\cdot)$ are given nonlinear functions and $D_c(\cdot)$ is a given matrix function of time.

The noise associated with system \eqref{eq:ct_x} - \eqref{eq:ct_y} can be described in terms of an IQC as in \cite{james1998,petersen1999},
\begin{equation}
\label{eq:iqc}
	\|x(0)-\bar{x}_0\|_N^2 + \int_{0}^{s} \left( \|w(t)\|_{Q_c}^2 + \|v(t)\|_{R_c}^2 \right) dt \leq d.
\end{equation}
Here, vectors are organized as columns,
$
    \|v\|_M:= \sqrt{\[v\]_M}
$
denotes the Euclidean (semi-) norm of a real vector $v$ generated by a real positive (semi-) definite symmetric matrix $M$, with
\begin{equation}
\label{quadro}
    \[v\]_M
    :=
    \[
        v, v
    \]_M,
    \qquad
    \[
        u,
        v
    \]_M
    :=
    u^{\rT}Mv.
\end{equation}
Unlike the semi-norm $\|\cdot\|_M$,  the quadratic form $\[\cdot\]_M$ is well-defined for any real symmetric matrix $M$. Also, $x(0)$ is the initial state value and $\bar{x}_0$ is the nominal initial state. A finite difference $(x(0)-\bar{x}_0)$ is allowed by a non-zero value of the constant $d$. If $d = 0$, then $x(0) = \bar{x}_0$. Also, $w(\cdot)$ and $v(\cdot)$ represent admissible uncertainties and $N = N^T > 0$ is a given matrix, $\bar{x}_0\in\mathbb{R}^n$ is a given state vector, $d>0$ is a given constant and $Q_c(\cdot)$, $R_c(\cdot)$ are given positive-definite, symmetric matrix functions of time.

In order to derive equations for a discrete-time set-valued state estimator, the continuous-time system in \eqref{eq:ct_x} - \eqref{eq:ct_y} needs to be discretized in reverse time. The reverse-time system formulation is used to formulate and solve the filtering problem which is recast as a subsidiary optimal control problem using the HJB equation. For further details, see \cite{ieeetac2009,acc2009}. Such a discretization can be achieved by using standard techniques such as the Euler or higher-order Runge-Kutta methods \cite{lewis2008}. In particular, applying the Euler scheme to \eqref{eq:ct_x} in reverse time yields,
\begin{align}
	x(t) \approx  x(t+\tau) - \tau a_c(t+\tau, x(t+\tau), u(t+\tau)) \nonumber\\
\label{eq:euler}
	\qquad  \qquad \qquad- \tau D_c(t+\tau) w(t+\tau),
\end{align}
where $\tau$ is the sampling time. Thus, \eqref{eq:euler} leads to the reverse-time discrete system of the form
\begin{align}
\label{eq:rt-x}
	x_k &= A_k(x_{k+1}) + B_k \: w_{k+1}, \\
\label{eq:rt-y}
	y_{k+1} &= C_k(x_{k+1}) + v_{k+1},
\end{align}
where $A_{(\cdot)}(\cdot)$ and $C_{(\cdot)}(\cdot)$ represent discrete-time nonlinear functions and $B_{(\cdot)}$ is a given time-varying matrix. The control variable $u$ is a known quantity and will be omitted for brevity as this paper deals with only the filtering problem.

Finally, the IQC in \eqref{eq:iqc} is discretized to obtain an equivalent SQC of the form
\begin{equation}
\label{eq:sqc}
	\|x_0-\bar{x}_0\|_N^2 + \sum_{s=0}^{T-1} \left( \|w_s\|_Q^2 + \|v_s\|_R^2 \right) \le d.
\end{equation}

\section{Set-valued State Estimation and the Optimal Control Problem}
\label{sec:svse}
Consider $y_k^0 = y_k$ to be a fixed measured output for the system \eqref{eq:rt-x} - \eqref{eq:rt-y} with disturbances bounded by the SQC \eqref{eq:sqc}. The set valued state estimation problem consists of constructing the set $Z_T[\bar{x}_0, y_{(\cdot)}^0 |_1^T,d]$ of all states $x_T$ at time step $T$ for the system \eqref{eq:rt-x} - \eqref{eq:rt-y} with initial conditions and disturbances defined by the quadratic constraint in \eqref{eq:sqc}, consistent with the measurement sequence $y_{(\cdot)}^0$.

Given an output sequence $y_{(\cdot)}^0$, it follows from the definition of $Z_T[\bar{x}_0, y_{(\cdot)}^0 |_1^T,d]$, that
\begin{equation}
	x_T \in Z_T[\bar{x}_0, y_{(\cdot)}^0 |_1^T,d]
\end{equation}
if and only if there exists a disturbance sequence $w_{(\cdot)}$ such that $J_T(x_T,w_{(\cdot)}) \le d$, where the cost functional $J_T(x_T,w_{(\cdot)})$ is obtained from the SQC \eqref{eq:sqc} and is of the form
\begin{align}
	J_T(\tilde x, w_{(\cdot)}) & \triangleq
	\frac{1}{2} \|{x}_0-\bar{x}_0\|_N^2 + \\ &\frac{1}{2}  \sum_{k=0}^{T-1}( \|w_k\|_{Q_k}^2  + \|v_{k+1}\|_{R_{k+1}}^2)
        \le d \label{eq:cost-functional}
\end{align}
with $v_{k+1} = y_{k+1}^0 - C_k(x_{k+1})$. Here, the vector $x_{(\cdot)}$ is the solution to the system \eqref{eq:rt-x} - \eqref{eq:rt-y} with input disturbance $w_{(\cdot)}$ and terminal condition $x_T = \tilde{x}$. Hence, 
\begin{equation}
\label{eq:zt}
	Z_T[\bar{x}_0, y_{(\cdot)}^0 |_1^T,d] = \left\{ \tilde{x} \in \mathbb{R}^{n} \,:\, \inf_{w_{(\cdot)}} J_T(\tilde{x}, w_{(\cdot)})\;\le \; d \right\}
	\hskip -1pt .
\end{equation}
The nonlinear optimal control problem for the system in \eqref{eq:rt-x} - \eqref{eq:rt-y} is defined by the optimization problem
\begin{align}
\label{eq:inf}
	V_T(\tilde x)\triangleq\inf_{w_{(\cdot)}} J_T(\tilde{x}, w_{(\cdot)}).
\end{align}
Here, it is assumed that the infimum in \eqref{eq:inf} exists. If not, the fulfillment of the inequality $J_T(x, w_{(\cdot)}) \le d$ does not guarantee the reachability of the terminal state $\tilde{x}$ under the SQC \eqref{eq:sqc}, in which case the inequality in \eqref{eq:zt} can only be defined as an inclusion,
\begin{equation}
\label{eq:zt2}
	Z_T[\bar{x}_0, y_{(\cdot)}^0 |_1^T,d] \subset \left\{ \tilde{x} \in \mathbb{R}^{n} \,:\, \inf_{w_{(\cdot)}} J_T(\tilde{x}, w_{(\cdot)})\;\le \; d \right\}.
\end{equation}

Now in order to obtain the optimal state estimates we  must  solve for the value function  \eqref{eq:inf}. This is done by applying the dynamic programming approach from optimal control theory.  In a discretized form, the value function satisfies the dynamic programming equation 
\begin{align}
V_{k+1}(x) = \min_{w_0} \Big\{ V_{k}(x(k-1)| x(k) = x))  \noncr
+  \frac{1}{2} \qform{w_0}{\qweight}  + \frac{1}{2}\|y - C(x)\|^2_R \Big\},
\end{align}
where $ V_{k}(x(k-1)| x(k) = x)) $ denotes the value function at time $k-1$ given a state $x$ at time $k$. 
Using the notation  $\oplus$ and $\otimes$ for the min-plus addition (min)  and multiplication (plus) 
operators respectively, we may rewrite the above as
\begin{align}
V_{k+1}(x) = \minplus_{w_0} \Big\{ V_{k}(x(k-1)| x(k) = x)  \noncr
\otimes \frac{1}{2} \qform{w_0}{\qweight}  \otimes \frac{1}{2}\|y - C(x)\|^2_R \Big\}.
\end{align}
In the following section we describe an approach  to solving  the above.
\section{Min-Plus structure preservation and filter design}
\label{sec:filter-equations}
In this section we solve the dynamic programming equation as follows. We express the value function in a particular min-plus basis (specifically the min of quadratic forms). Then we exploit the linearity of the dynamic programming operator in this space to obtain a recursive equation for the parameters used in this expansion.   This recursion is possible owing to the fact that after  propagation by the dynamic programming operator, this min-of-quadratic-forms structure is preserved. The submatrices of the quadratic form, in fact, correspond to the solution of the  Riccati  equation for optimal filter design.


%
We will omit the time subscripts for the state and nonlinear functions for brevity.

From \eqref{eq:cost-functional} and \eqref{eq:inf} we have at $T=0$ 
\begin{align}
\label{eq:v0}
V_0(x) := \frac{1}{2}  \Bigg\{ \|x - \bar{x}_0\|_{N^0}^2 + \phi^0 \Bigg\},
\end{align}
which can be written in the quadratic form
\begin{equation}
\label{eq:v0-quad}
	\bigwedge_{i \in \cI_0}\frac{1}{2}  \qformvec{x}{Q^{v,0}_i},
\end{equation}
where
\begin{equation}
\label{\eqnname defofQ}
Q^{v,0}_i 
:= \left[\begin{array}{cc}N^0_i & {L^0_i}^T \\{L^0_i} & \bar{\phi}^0_i\end{array}\right],
\end{equation}
$L_1^0= -  {\bar{x}_0}^T N^0$, $\bar\phi^0_1=\bar x_0^T N^0\bar x_0+\phi^0$ and $\cI_0=\{1\}$.

At $T=1$, the dynamic programming recursion equation can be written in the form
\begin{align}
V_1(x) := \minplus_{w_0} \Big\{V_0&(A(x) + B w_0) + \frac{1}{2} \|w_0\|_{\qweight}^2 \nonumber \\
	&+ \frac{1}{2} \|y-C(x)\|_R^2  \Big\}, 
\label{\eqnname dpprecursion}
\end{align}
where we let $A$ and $B$ be time-independent for simplicity.

Substituting for $V_0$ from \eqref{eq:v0} - \eqref{\eqnname defofQ} in \eqref{\eqnname dpprecursion} and using the backward time dynamics \eqref{eq:rt-x} we obtain,
\begin{align}
	&V_1(x) \nonumber \\
	&= \minplus_{w_0} \left(\bigwedge_{i \in \cI_0}  \frac{1}{2}  \|f_{x,w_0}\|_{Q^{v,0}_i}^2 
	+  \frac{1}{2}  \|w_0\|_{\qweight}^2 + \frac{1}{2} \|g_{y,x}\|_R^2 \right)\hskip -2pt, \nonumber \\
	&\hskip -1pt=\hskip -3pt  \bigwedge_{i \in \cI_0}\hskip -3pt \left\{ \minplus_{w_0}  \frac{1}{2}  \left( \|f_{x,w_0}\|_{Q^{v,0}_i}^2 
	+ \|w_0\|_{\qweight}^2 \right)\hskip -1pt + \frac{1}{2} \|g_{y,x}\|_R^2 \right\}\hskip -2pt, \noncr
\label{\eqnname dpp2}
\end{align}
where, for the sake of enhancing clarity,  we use  the notation 
\begin{align}
f_{x,w_0}^T &:= \left[ (A(x) + B w_0)^T \;\;\; 1 \right],\\
g_{y,x}&:= y-C(x).
\end{align}

The minimizing $w_0$ is found from the following expression:
\begin{align}
\argmin_{w_0}\left[ \| f_{x,w_0} \|_{Q^{v,0}_i}^2 + \|w_0\|_{Q_{\eta}}^2 \right].
\end{align}
Solving for $w_0$ and rewriting the matrix $Q^{v,0}_i$ in terms of its constituent matrices from \eqref{\eqnname defofQ} yields
\begin{align}
w^*_0 &= - [\qweight +  B^T N^0_i B]^{-1}\times [B^T {{L^0_i}^T} +  B^T N^0_i A(x)]
\label{eq:w-opt}
\end{align}
which is of the form $ w^i_c + w^i_l A(x)$. Substituting $w_0 = w_0^*$ from \eqref{eq:w-opt} in \eqref{\eqnname dpp2} we obtain
\begin{align}
V_1(x)  &= \bigwedge_{i \in \cI_0}  \frac{1}{2}  \|\tilde{f}_{x,w_0^*}\|_{Q^{v,0}_i}^2 +  \frac{1}{2}  \| {w_0^*} \|_{\qweight} + \frac{1}{2} \|g_{y,x}\|^2_R,
\end{align}
where $\tilde{f}_{x,w_0}^T = \left[\left(A(x)+B w^i_c + B w^i_l A(x) \right)^T \;\;\; 1 \right]$ and $\tilde{w} = [w^i_c + w^i_l A(x)]$. Collecting terms in $A(x)$ we find
\begin{align}
	&V_1(x) \nonumber \\
	&= \bigwedge_{i \in \cI_0} \Bigg\{  \frac{1}{2}  \|A(x)\|_{h_{w_l^i}}^2 + \Big\{L^0_i (I + B w^i_l) \nonumber \\
	&+  {w^i_c}^T B^T N^0_i (I + B w^i_l) +  {w^i_c}^T \qweight w^i_l \Big\} A(x)  \nonumber \\
	&+ L^0_i B w^i_c + \frac{1}{2} \Bigg[[{w^i_c}^T B^T N^0_i B w^i_c] + \qform{(w^i_c)}{\qweight} \noncr & \quad \quad\quad \quad\quad\quad\quad+ \bar{\phi}^0_i  \Bigg] \Bigg\}
	+ \frac{1}{2}{\|y - C(x)\|}_R^2,
\label{\eqnname v1expandedform}
\end{align}
where $h_{w_l^i} := \qform{(I + B w^i_l)}{N^0_i} + \qform{w^i_l}{\qweight}$.
Consider the following quadratic approximations
\begin{align}
- \braket{y} {C(x)}_R &= \bigwedge_{j \in \cJ} \qformvec{x} {|y| Q^{c,y}_j} \label{eq:y}\\
 \|C(x)\|^2_R &:= \bigwedge_{ l \in \cL} \qformvec{x}{Q^b_l}.\label{eq:cx}
\end{align}
Note the dependence of the above terms on the output $y$ (however more specifically $Q^{c,y}_j$ only depends on 
the sign of $y$). 
Adding \eqref{eq:y} and \eqref{eq:cx} yields
\begin{align}
- \braket{y}{C(x)}_R &+ \|C(x)\|^2_R \nonumber \\
	&= \bigwedge_{m \in \cM} \qformvec{x}{Q^0_m},
\end{align}
where $Q^0_m$ is a matrix, some of whose terms depend on $y$. Here  $\cM := \cJ \times \cL$ for which the following holds: $\forall j \in \cJ, \,l \in \cL, \,\,\,\exists m \in \cM$ such that 
\begin{align}
 Q^0_m :=    |y|\,Q^{c,y}_j + Q^b_l.
\end{align}
We further have the following representations for the terms in \eqref{\eqnname v1expandedform}:
\begin{align}
	&A(x)^T M^0_i A(x) := \bigwedge_{a \in \cI^0_a} \qformvec{x}{Q^a_0(M^0_i))}, \noncr
	&\tilde{M}^0_i A(x) := \bigwedge_{\tilde{a} \in {\tilde{\cI}}^0_a} \qformvec{x}{Q^{\tilde{a}}_0(\tilde{M}^0_i))}.\label{eq:qc}
\end{align}
Thus using \eqref{eq:y} - \eqref{eq:qc} in \eqref{\eqnname v1expandedform} for $V_1(x)$ we obtain
\begin{align}
V_1(x) 	&= \bigwedge_{i \in \cI_0}\Bigg\{\bigwedge_{a \in \cI^0_a} \qformvec{x}{Q^a_0(M^0_i)} \nonumber \\
	&+   \bigwedge_{\tilde{a} \in {\tilde{\cI}}^0_a} \qformvec{x}{Q^{\tilde{a}}_0(\tilde{M}^0_i))} \nonumber \\
	&+ \qformvec{x}{ Q^c } \nonumber \\
	&+ \bigwedge_{m \in \cM} \qformvec{x}{Q^0_m} \Bigg\}, \label{\eqnname v1longand}
\end{align}
where
\begin{align}
\phi^1_i &:= \bar{\phi}^0_i + 2 L^0_i B w^i_c + [{w^i_c}^T B^T N^0_i B w^i_c] \nonumber \\  
&\qquad \qquad\qquad + \qform{(w^i_c)}{\qweight},  \nonumber \\
	Q^c &:= \left(\begin{array}{ccc}0 & 0 & 0 \\0 & 0 &0 \\0 & 0 & \phi^1_i\end{array}\right).
\end{align}
Now, to further simplify the  expression \eqref{\eqnname v1longand}, we note that for all $i \in \cI_0, \, a \in \cI^0_a,\, \tilde{a} \in \tilde{\cI}^0_a, \, m \in \cM, \, \exists k \in \cI_1$  (where $\cI_1:=  \cI_0 \times \cI^0_a \times \tilde{\cI}^0_a \times  \cM$) such that
\begin{align}
Q^{v,1}_k := Q^a_0(M^0_i) + \tilde{Q}^{\tilde{a}}_0(\tilde{M}^0_i)  + Q^c + Q^0_m.
\end{align}

Hence \eqref{\eqnname v1longand} can be written in the form
\begin{align}
 \bigwedge_{k \in \cI_1} \qformvec{x}{Q^{v,1}_k}.
\end{align}
Thus we have a recursive relationship in the \textit{coefficients of the quadratic forms}  between two consecutive time steps. By propagating these terms 
across multiple time steps  we may evaluate the cost function at any desired $x$ (without storing the value for $x$ at each time step).

From the results above, after performing the minimization with respect to the set of quadratics at the current time step $t$, the value function  has the form
\begin{align}
V_t(x) = \bigwedge_{k \in \cI_N} \frac{1}{2}  \qformvec{x}{Q^{v,N}_k}.
\label{eq:vt-minplus}
\end{align}
To obtain a state estimate $\hat{x}^*$ we minimize the above with respect to $x$, i.e.,
\begin{align}
\hat{x}^* = \argmin_{x} V_t(x). \label{\eqnname xhatstardef}
\end{align}
A set-valued estimate is obtained as a sub-level set of $V_t(\cdot)$.
Note that the minimizing $\hat{x}^*$ occurs  at one of the troughs of one of the quadratics (say $Q^*$) which has the structure
\begin{align}
 \left[\begin{array}{cc}N^*_i & {L^*}^T \\{L^*} & \bar{\phi}^*\end{array}\right].
\end{align}
As indicated in \cite{ieeetac2009,acc2009} the value function may also be associated with the real symmetric precision matrix $\Pi_t \succ 0$ and the state estimate $\hat{x}^*_t$ as follows
\begin{align}
V_t(x) = \frac{1}{2}  \qform{(x  - \hat{x}^*_t)} {\Pi_t} + \frac{1}{2}  \hat{\phi}_t.
\label{eq:rekf}
\end{align}
By comparing coefficients in \eqref{eq:vt-minplus} and \eqref{eq:rekf} it can be seen that
\begin{align}
\Pi_t &= N^*,\noncr
\hat{x}^*_t &= -\frac{1}{2} [{N^*}^T]^{-1} \,{L^*}^T.
\end{align}
Here, the matrix $\Pi^{-1}_t = P_t$ corresponds to the estimation error covariance matrix in the traditional Kalman filter in the stochastic setting.

Note that in order to choose the minimizing quadratic we obtain the minimizing point $x^*$ for each quadratic as follows. Given a form
\begin{align}
\frac{1}{2} [x^T\,\, 1] \left(\begin{array}{cc}q_{11} & q_{12} \\q_{21} & q_{22}\end{array}\right)  \left[\begin{array}{c}x \\1\end{array}\right],
\end{align}
the minimizing $x^*$ for this quadratic is given by 
\begin{align}
x^* = -[q_{11} + q_{11}^T]^{-1} [q_{12} + q_{21}^{T}].
\label{eq:xstar}
\end{align}
The latter is true if and only if the states are free to take on any values. In the case where the states are constrained, the minimization in \eqref{\eqnname xhatstardef} must be performed in the permissible set of states.

\section{Illustrative Example}
\label{sec:ex}
In order to demonstrate  the concepts introduced in this article, we analyze a two dimensional system with linear dynamics and a nonlinear output function, defined as follows
\begin{align}
\frac{d}{dt} \left[\begin{array}{c}x_1 (t) \\x_2 (t) \end{array}\right]  &= \left[\begin{array}{cc}0 & 0 \\1 & 0\end{array}\right]  \left[\begin{array}{c}x_1 (t) \\x_2 (t) \end{array}\right]  + \left[\begin{array}{c}1 \\0\end{array}\right] w(t),\\
y(t) & = \sqrt{2} \sin(x_2 (t))  + v(t),
\end{align}
where $w(\cdot)$ and $v(\cdot)$ are the process disturbance and measurement noise respectively.
After discretization the system dynamics is
\begin{align}
\left[\begin{array}{c}x_1 (k+1) \\x_2 (k+1) \end{array}\right]  &= \left[\begin{array}{cc}1 & 0 \\0.1 & 1\end{array}\right]  \left[\begin{array}{c}x_1 (k) \\x_2 (k) \end{array}\right]  \\ & \quad + \left[\begin{array}{c}0.1 \\0\end{array}\right] \Delta B_k, \label{\eqnname discdynamicsexample}
\end{align}    
where $\Delta B_k$ is the increment corresponding to $w(\cdot)$ over the sampling time.
    
\begin{figure}[htp]
   \centering
   \subfloat[Approximating $+ \sin(\theta)$]{\includegraphics[width= .7\hsize]{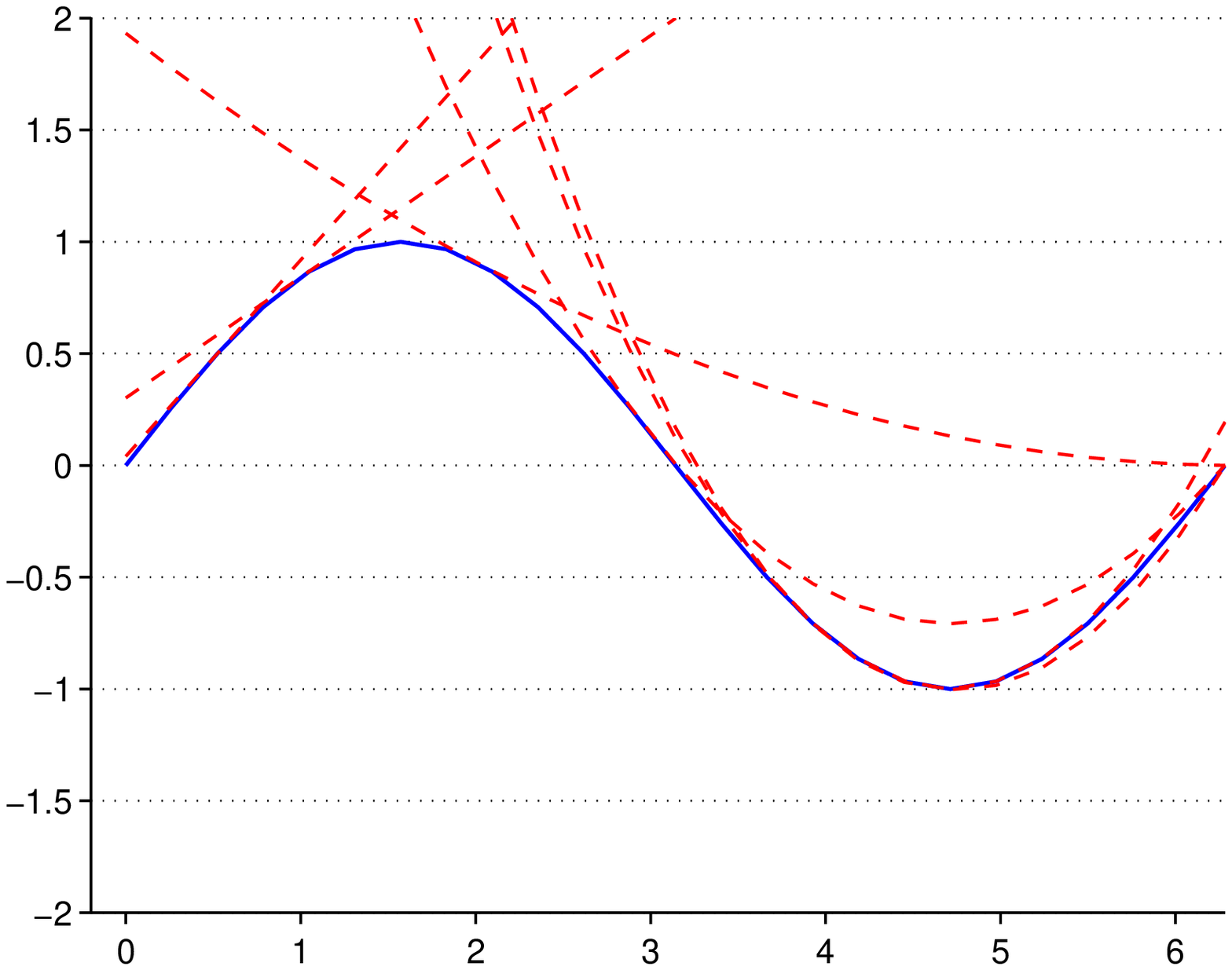}}       \\   
  \subfloat[Approximating $- \sin(\theta)$]{\includegraphics[width=.7 \hsize]{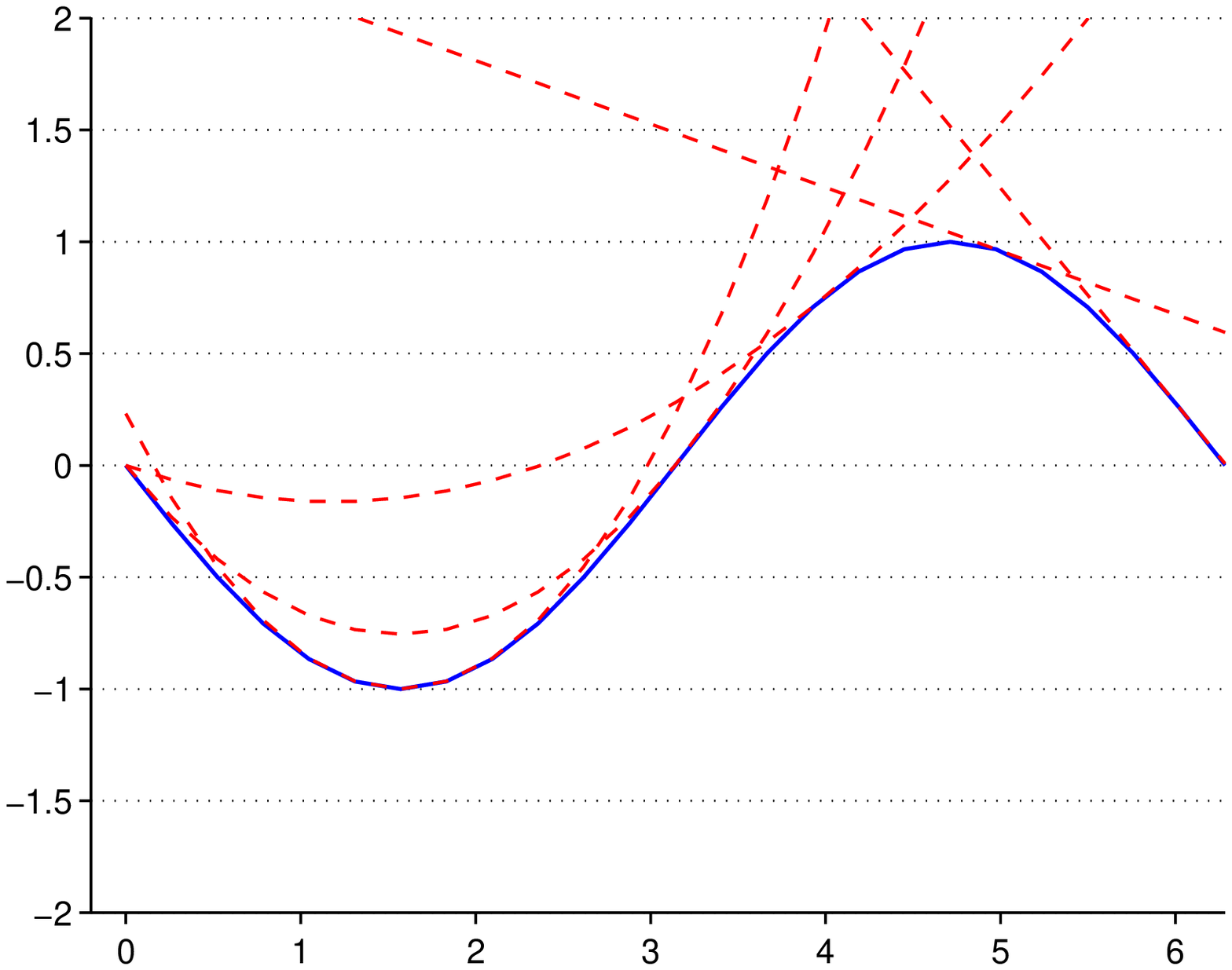}}
  \caption{Approximation of $\pm \sin(\theta)$ as a min of quadratics.}
    \label{\figname figcapprox}
\end{figure}

\begin{figure}[htp]
   \centering
  \includegraphics[width= .8\hsize]{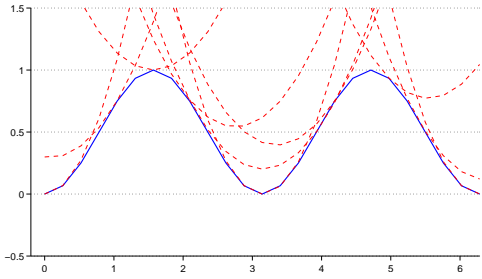}
    \caption{Approximation of $\sin^2(\theta)$ as a min of quadratics.}
     \label{\figname figccsqapprox}
\end{figure}

In order to apply the deterministic filtering approach we approximate the output function $C(\theta) := \pm \sqrt{2} \sin(\theta)$ (where the  sign of the function used depends on the sign of the output as described previously in \eqref{eq:y}  and $[C(\theta)]^2$ as a minimum of convex functions as indicated in Fig.~\ref{\figname figcapprox}, \ref{\figname figccsqapprox}.

%

\begin{figure}
\begin{center}
\includegraphics[width= \hsize]{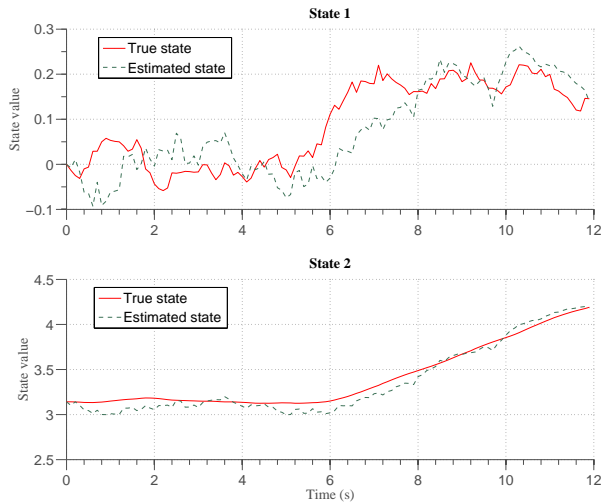}
\caption{State filtering}
\label{\figname figstateest1}
\end{center}
\end{figure}

\begin{figure}
\begin{center}
\includegraphics[width= \hsize]{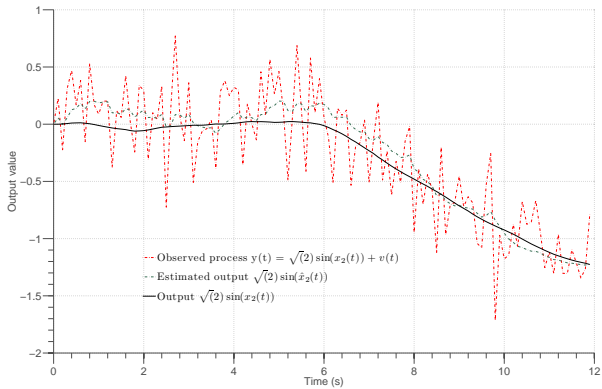}
\caption{Filtered measurement}
\label{\figname figoutputest1}
\end{center}
\end{figure}

By applying the min-plus filter design  approach we obtain the estimation results indicated in Fig.~\ref{\figname figstateest1}, \ref{\figname figoutputest1}.   Intuitively, the first state is more difficult to estimate  -- as can be inferred from \eqref{\eqnname discdynamicsexample}, there is a weak dependence of the second state on the first (in addition to a one sample delay) and the noise increment in the output has a reasonably high variance. 

%

\section{Conclusion and Future Directions}
\label{sec:conclusions}
The technique described herein provides an approach to the design of filters for systems with nonlinear dynamics and nonlinear output. Its main contribution is in the utilization  of the min-plus basis expansion of the value function coupled with the exploitation of the linearity of the dynamic programming operator over such a (semi)-field. A few  of the avenues along which a study of the  ramifications and salient features of these methods may be pursued are: the error analysis of the dependence of the accuracy on the approximation of the output and system dynamics  by convex functions,  the extension to systems with uncertainty, and the development of optimal approximation techniques for approximation of any desired function via a sequence of convex functions.  In addition, these methods provide a computationally tractable approach for  nonlinear filtering, and the applications of this to time critical problems would also provide a fruitful direction of practical relevance while driving further insights into these classes of approaches.

\bibliographystyle{unsrt}
\balance

%

\end{document}